\documentclass[12pt,twoside]{article}
\usepackage{latexsym}
\usepackage{mathrsfs}
\usepackage{graphicx}
\usepackage{arydshln}
\def\thebibliography#1{\center{\bf\normalsize References}\list
 {[\arabic{enumi}]}{\settowidth\labelwidth{[#1]}\leftmargin\labelwidth
 \advance\leftmargin\labelsep
 \usecounter{enumi}}
 \def\newblock{\hskip .11em plus .33em minus .07em}
 \sloppy\clubpenalty4000\widowpenalty4000
 \sfcode`\.=1000\relax}
\def\tt{\hspace{-0.19cm}{\bf .}\hspace{0.19cm}}

\def\sh{\,{\rm sh}\,}

\def\proof{\bb{\bf Proof.\hspace{0.19cm}}}
\def\endpf{\hfill$\Box$\vspace{0.4cm}}
\def\nnb{\nonumber}
\def\disp{\displaystyle}
\newcommand{\eqref}[1]{$(\ref{#1})$}
\newcommand{\refeq}[1]{$(\ref{#1})$}

\def\Ga{\alpha}
\def\Gb{\beta}
\def\Gg{\gamma}
\def\Gd{\delta}

\def\Gl{\lambda}

\def\GL{\Lambda}

\def\bb{\bar b}

\def\mm{\mbox{ }}

\def\bb{\hspace{18pt}}

\def\qq{\qquad}
\def\q{\quad}

\def\eqif{\mm {\rm if}\mm }

\def\eqin{ \mm {\rm in } \mm }

\def\all{  \mm \forall \mm }

\newcommand{\dpp}[2]{{\partial {#1} \over \partial {#2}}}
\newcommand{\doo}[2]{{d {#1} \over d {#2}}}
\newcommand{\ddp}[3]{{\partial {^{#1}#2} \over \partial {#3^{#1}}}}

\def\sh{\mm{\rm sh}\mm}

\def\defeq{\stackrel{\triangle}{=}}

\begin{document}

\title{
{\bf  The Proof of Alzer's Conjecture on Generalized Logarithmic
Mean }
\thanks{The authors were supported in part by NSFC (No. 10671040) and
FANEDD (No. 200522).} }
\author{ Hongwei Lou\thanks{School of Mathematical Sciences,
Fudan University, Shanghai 200433, China, \& Key Laboratory of
Mathematics for Nonlinear Sciences (Fudan University), Ministry of
Education   (hwlou@fudan.edu.cn). } and Dongdi Liu\thanks{School of Mathematical Sciences,
Fudan University, Shanghai 200433, China (0418127@fudan.edu.cn). } }
\date{}
\maketitle
\begin{quote}
\footnotesize {\bf Abstract.} In 1987, Alzer posed a
conjecture on generalized logarithmic mean, which was introduced by
Stolarsky in 1975. To prove Alzer's conjecture,
Lou posed a conjecture on generalized
inverse harmonic mean in 1995.
By proving Lou's conjecture,  the paper yields
  Alzer's  conjecture finally.



\end{quote}
\normalsize
\newtheorem{Definition}{\bb Definition}[section]
\newtheorem{Theorem}[Definition]{\bb Theorem}
\newtheorem{Lemma}[Definition]{\bb Lemma}
\newtheorem{Corollary}[Definition]{\bb Corollary}
\newtheorem{Proposition}[Definition]{\bb Proposition}
\newtheorem{Remark}{\bb Remark}[section]
\newtheorem{Example}{\bb Example}[section]
\newfont{\Bbb}{msbm10 scaled\magstephalf}
\newfont{\frak}{eufm10 scaled\magstephalf}
\newfont{\sfr}{msbm7 scaled\magstephalf}

\def\theequation{1.\arabic{equation}}
\setcounter{equation}{0} \setcounter{section}{1}
\setcounter{Definition}{0}\setcounter{Example}{0} \textbf{1.
Introduction.} For two positive numbers $a$ and $b$, Stolarsky
defined in \cite{S} the generalized logarithmic mean of $a$, $b$ as
\begin{equation}\label{E101}
L_r(a,b)\defeq \Big({b^r-a^r\over r(b-a)}\Big)^{1\over r-1},
\end{equation}
where $r\in [-\infty,+\infty]$ and $L_{-\infty}(a,b)$, $L_0(a,b)$,
$L_1(a,b)$, $L_{+\infty}(a,b)$ are looked as the corresponding
limits:
$$
L_{-\infty}(a,b)\defeq \lim_{r\to -\infty}L_r(a,b)=\min (a,b),
$$
$$
L_0(a,b)\defeq \lim_{r\to 0}L_r(a,b)={b-a\over \ln b-\ln a},
$$
$$
L_1(a,b)\defeq \lim_{r\to 1}L_r(a,b)={1\over e}\Big({b^b\over
a^a}\Big)^{1\over b-a},
$$
$$
L_{+\infty}(a,b)\defeq \lim_{r\to +\infty}(a,b)=\max (a,b).
$$
Similarly, in this paper, the value of a function on its contact discontinuity point is always looked
as its corresponding limit.  The generalized logarithmic mean has been studied by many
researchers and it is still an interesting topic today (see
\cite{A2}---\cite{LP}, \cite{NP}---\cite{Q},
\cite{SY1}---\cite{SY2}, for examples). The aim of this paper is to
prove the following inequalities related to generalized logarithmic
mean:
\begin{equation}\label{E102}
2L_0(a,b)<L_r(a,b)+L_{-r}(a,b)<a+b,\q\all r\in (0,+\infty), b>a>0.
\end{equation}
The above inequalities is a conjecture posed by Alzer \cite{A} in
1987. Alzer himself proved that
\begin{equation}\label{E103AA}
L_1(a,b)+L_{-1}(a,b)>2L_0(a,b),\qq   \all b>a>0
\end{equation}
and the following result:
\begin{Lemma}\tt\label{T101} For any $r\in (0,+\infty)$, $b>a>0$, it holds that
\begin{equation}\label{E103A}
ab<L_r(a,b)L_{-r}(a,b)<L^2_0(a,b).
\end{equation}
\end{Lemma}

To prove \refeq{E102}, Lou studied generalized inverse harmonic
mean (which is a special case of Gini mean \cite{Gini}) of two
positive numbers in \cite{Lou},
\begin{equation}\label{E103}
C_r(a,b)\defeq \Big({b^r+a^r\over b+a}\Big)^{1\over r-1},
\end{equation}
where $r\in [-\infty,+\infty]$.
We mention that
$$
C_0(a,b)=L_2(a,b)={a+b\over 2}, \q C_{-1}(a,b)=L_{-1}(a,b)=\sqrt{ab}, \q
C_2(a,b)={a^2+b^2\over a+b}
$$
are the arithmetic mean, the geometric mean and the inverse harmonic
mean, respectively.  While
$$
C_{-\infty}(a,b)=\min (a,b),\q
C_1(a,b)=\Big({b^b a^a}\Big)^{1\over
b+a},\q
C_{+\infty}(a,b)=\max (a,b).
$$
On the other hand, we have
\begin{equation}\label{E104}
L_r(a^2,b^2)=L_r(a,b)C_r(a,b),\qq\all r\in [-\infty,+\infty], a,b>0.
\end{equation}
Using \refeq{E104}, Lou observed in \cite{Lou} that \refeq{E102} can
be proved if the following equalities hold (see the proofs of
Conjectures A and B in Sections 3 and 5):
\begin{equation}\label{E105}\left\{\begin{array}{l}\disp
C_r(a,b)+C_{-r}(a,b)>a+b, \\
\disp C_r^2(a,b)+C^2_{-r}(a,b)<a^2+b^2,\end{array}\right.\q\all
r\in (0,+\infty), b>a>0.
\end{equation}
More precisely, rewrite \refeq{E102} and \refeq{E105} as\\
\textbf{Conjecture A (H. Alzer)} It holds that
\begin{equation}\label{E106}
L_r(a,b)+L_{-r}(a,b)> 2L_0(a,b),\qq\all r\in (0,+\infty], b>a>0;
\end{equation}

\noindent
 \textbf{Conjecture B (H. Alzer)}  It holds that
\begin{equation}\label{E107}
L_r(a,b)+L_{-r}(a,b)< a+b,\qq\all r\in [0,+\infty), b>a>0;
\end{equation}

\noindent\textbf{Conjecture 1 (H. Lou)}  It holds that
\begin{equation}\label{E108}
C_r(a,b)+C_{-r}(a,b)>a+b,\qq\all r\in (0,+\infty), b>a>0;
\end{equation}
\noindent \textbf{Conjecture 2 (H. Lou)}  It holds that
\begin{equation}\label{E109}
C_r^2(a,b)+C^2_{-r}(a,b)<a^2+b^2,\qq\all r\in [0,+\infty), b>a>0.
\end{equation}
Then, noting that
\begin{equation}\label{E110}
\lim_{b\to a}{L_r(a,b)+L_{-r}(a,b)- 2L_0(a,b)\over
(b-a)^4}={r^2\over 960a^3}>0,
\end{equation}
and
\begin{equation}\label{E110A} \lim_{b\to
a}{L_r(a,b)+L_{-r}(a,b)- a-b\over (b-a)^2}=-{1\over 6a}<0,
\end{equation}
Lou showed in \cite{Lou} that Conjecture 1 implies Conjecture A while
Conjecture 2 implies Conjecture B.
Unfortunately, Conjectures 1 and
  2 are also difficult to prove though some special cases
were verified in \cite{Lou}. It was proved there that Conjectures 1
and  A hold when $r=1, 2,{1\over 2},3,{1\over 3},{3\over
2}, {2\over 3}$, while  Conjectures 2 and  B hold when
$r\in [{1\over 7},7]$.

In this paper, by the help of symbolic calculation in computer, we are able to prove Conjectures 1 and 2. And then
we get the proofs of Conjectures A and B.

We would like to mention that since the
Stolarsky mean (\cite{S})
$$
E_{p,q}(a,b)\defeq \Big({q\over p} {b^p-a^p\over
b^q-a^q}\Big)^{1\over p-q}
$$
and the Gini mean
$$
G_{p,q}(a,b)\defeq \Big( {b^p+a^p\over b^q+a^q}\Big)^{1\over p-q}
$$
can be got by
$$
E_{p,q}(a,b)=\left\{\begin{array}{ll}\disp
\Big(L_{p/q}(a^q,b^q)\Big)^{1\over q}, & \eqif q\ne 0,\\
\disp L_0(a,b), & \eqif p=q=0,
\end{array}\right.
$$
and
$$
G_{p,q}(a,b)= \left\{\begin{array}{ll}\disp
\Big(C_{p/q}(a^q,b^q)\Big)^{1\over q},, & \eqif q\ne 0,\\
\disp C_0(a,b), & \eqif p=q=0,
\end{array}\right.
$$
in some sense, it is enough to study $L_r(a,b)$ and $C_r(a,b)$ when
one need to study $E_{p,q}(a,b)$ and $G_{p,q}(a,b)$.

Sections 2 and 4 are devoted to prove Conjectures 1 and 2,
while Sections 3 and 5 are devoted to prove Conjectures A and B.

\if{
Lou proved Conjecture 1 and therefore proved Conjecture A. However, Lou failed to prove Conjecture 2.
To prove Conjecture B, Lou generalize $C_r(a,b)$ in the
following manner:
\begin{equation}\label{E502}
H_{r,s}(a,b)\defeq {L_r(a^{s+1},b^{s+1})\over
L_r(a^s,b^s)}=\Big({(b^{(s+1)r}-a^{(s+1)r})(b^s-a^s)\over
(b^{s+1}-a^{s+1})(b^{sr}-a^{sr})}\Big)^{1\over r-1}.
\end{equation}
Then, Lou proved  the following proposition and then get a proof of Conjecture B:
\begin{Proposition}\tt\label{T502}
Let $s\geq 1$, $\Gd\geq 1$,  $a,b>0$, $a\ne b $, $0\leq r<+\infty$,
$\Ga>0$, Then, as a function of $r$,
$$
 \Big(\max
(a,b)\Big)^{s+\Gd}H^\Ga_{r,s}(a,b)+\Big(\min
(a,b)\Big)^{s+\Gd}H^\Ga_{-r,s}(a,b)
$$
is increasing strictly in $[0,+\infty)$. In particular,
\begin{equation}\label{E506}
 \Big(\max
(a,b)\Big)^{s+\Gd}H^\Ga_{r,s}(a,b)+\Big(\min
(a,b)\Big)^{s+\Gd}H^\Ga_{-r,s}(a,b)<a^{s+\Gd+\Ga}+b^{s+\Gd+\Ga}.
\end{equation}
\end{Proposition}
Later, Liu gave  a proof of Conjecture 2 in his  thesis and then got a simple proof of
Conjecture B based on Conjecture 2.

}\fi

\def\theequation{2.\arabic{equation}}
\setcounter{equation}{0} \setcounter{section}{2}
\setcounter{Definition}{0} \setcounter{Example}{0}\textbf{2.  Proof
of Conjecture 1.}  We recall some basic properties of $L_r(a,b)$  and $C_r(a,b)$.
\begin{Proposition}\tt\label{T200} Assume $a,b>0$, $r\in [-\infty,+\infty]$.
\begin{enumerate}\renewcommand{\labelenumi}{{\rm (\roman{enumi})}}
\item $L_r(a,b)$ is symmetric, that is,
\begin{equation}\label{E201}
L_r(a,b)=L_r(b,a).
\end{equation}
\item For any $\Ga>0$,
\begin{equation}\label{E202}
L_r(\Ga a,\Ga b)=\Ga L_r(a,b).
\end{equation}
\item For any $-\infty< s<r<+\infty$, $b>a>0$,
\begin{equation}\label{E303A}
\min (a,b)<L_s(a,b)<L_r(a,b)<\max (a,b).
\end{equation}
\end{enumerate}
\end{Proposition}
The proof of the above proposition can be found in \cite{S}.

\begin{Proposition}\tt\label{T301} Assume  $a,b>0$, $r\in [-\infty,+\infty]$.
\begin{enumerate}\renewcommand{\labelenumi}{{\rm (\roman{enumi})}}
\item $C_r(a,b)$ is symmetric, that is,
\begin{equation}\label{E301}
C_r(a,b)=C_r(b,a).
\end{equation}
\item For any $\Ga>0$,
\begin{equation}\label{E302}
C_r(\Ga a,\Ga b)=\Ga C_r(a,b).
\end{equation}
\item For any $-\infty< s<r<+\infty$, $b>a>0$,
\begin{equation}\label{E303}
\min (a,b)<C_s(a,b)<C_r(a,b)<\max (a,b).
\end{equation}
\item Let  $0<r<s<+\infty$, $b>a>0$. Then
\begin{equation}\label{E304}
C_{-1}^2(a,b)<C_s(a,b)C_{-s}(a,b)<C_r(a,b)C_{-r}(a,b)<C^2_0(a,b).
\end{equation}
\end{enumerate}
\end{Proposition}
\proof Though the proof of above proposition was given in
\cite{Lou}, for the convenience of readers, we give the proofs of
(iii)---(iv) in the following. Without loss of generality, we set
$b>a=1$.

\noindent (iii) It suffice to prove that
$$
f(r,b)\defeq \dpp {} r \Big(\ln C_r(1,b)\Big)=-{1\over (r-1)^2}\ln
{b^r+1\over b+1}+{1\over r-1}{b^r\ln b\over b^r+1}.
$$
is positive. Denote
$$
g(r,b)=(1-r)^2f(r,b).
$$
We have
$$
\dpp {g(r,b)} r=(r-1){b^r\ln^2 b\over (b^r+1)^2}.
$$
Thus, for fixed $b>1$, $g(r,b)$ is decreasing strictly in $r\in (0,1)$ and increasing strictly in $r\in (1,+\infty)$.
Therefore,
$$
g(r,b)>g(1,b)=0, \qq\all r\ne 1.
$$
Consequently, $f(r,b)$ is positive since
\begin{eqnarray*}\disp
 f(1,b)&=& \lim_{r\to 1}f(r,b) = {b\ln^2 b\over 2(b+1)^2}>0.
\end{eqnarray*}
(iv) We have
$$
\dpp{} r \Big[\ln \Big(C_r(1,b)C_{-r}(1,b)\Big)\Big]=f(r,b)-f(-r,b).
$$
Let
\begin{eqnarray*}
h(r,b)&=& {(r-1)^2(r+1)^2\over r}(f(r,b)-f(-r,b))\\
&=& {(r+1)^2\over
r}g(r,b)-{(r-1)^2\over r}g(-r,b), \qq r>0.
\end{eqnarray*}
Then we can get that
\begin{eqnarray*}\disp
 \dpp {h(r,b)} r &= &{(r^2-1)\ln b\over r^2
(b^r+1)^2}\Big(1-b^{2r}+2rb^r\ln b\Big) \\
\disp &=& {(1-r^2) \ln^2 b\over r (b^r+1)^2}\Big(L_0(1,b^{2r})-L_{-1}(1,b^{2r})\Big).
\end{eqnarray*}
\if{
\begin{eqnarray*}\disp
 \doo {h(r)} r &= &  {r^2-1\over r^2} \Big(g(r)-g(-r)\Big)+ {(r+1)^2\over
r}(r-1){b^r\ln^2b\over (b^r+1)^2}\\
\disp & & -{(r-1)^2\over r}(r+1){b^r\ln^2b\over (b^r+1)^2}\\
&=& \Big(1-{1\over r^2}\Big){ \ln b\over
(b^r+1)^2}(1-b^{2r}+2b^r\ln b^r) \\
\disp &=& -2\Big(1-{1\over r^2}\Big){  \ln b^r\ln b\over
(b^r+1)^2}\Big(L_0(1,b^{2r})-L_{-1}(1,b^{2r})\Big).
\end{eqnarray*}}\fi
Thus, $h(r,b)$ is increasing strictly in $r\in (0,1)$ and decreasing strictly in $r\in (1,+\infty)$.
Consequently,
$$
h(r,b)<h(1,b)=0,  \qq\all r>0, r\ne 1.
$$
Therefore, $ \ln \Big(C_r(1,b)C_{-r}(1,b)\Big)$ is decreasing strictly in $r\in (0,+\infty)$ and \refeq{E304} follows.
\endpf

Now, we begin to prove Conjecture 1 and state it as
\begin{Theorem}\tt\label{T303} Let $r\in (0,+\infty)$, $b>a>0$. Then
we have
\begin{equation}\label{E305}
C_r(a,b)+C_{-r}(a,b)>a+b.
\end{equation}
\end{Theorem}
\proof Without loss of generality, we set $b>a=1$. Since
$$
\Big(C_{{1\over r}}(1,x)+C_{-{1\over
r}}(1,x)-x-1\Big)\Big|_{x=b^r}=\Big(C_r(1,b)+C_{-r}(1,b)-b-1\Big){b^r+1\over
b+1},
$$
we see that \refeq{E305} holds for some $r=r_0\in (0,1]$ if and only
if it holds for $r={1\over r_0}$. Therefore, we can suppose that
$r\geq 1$ without loss of generality.
We have
\begin{eqnarray}\disp \label{E306}
 \dpp {C_r(1,b)} b  &=& 
{rb^r+rb^{r-1}-b^r-1\over (r-1)(b^r+1)(b+1)}C_r(1,b),
\\
\label{E307}\disp \dpp {C_{-r}(1,b)} b &=&  {r +rb^{ -1}+b^r+1 \over
 (r+1)(b^r+1)(b+1)}C_{-r}(1,b).
\end{eqnarray}
Let
$$
F_0(r,b)={C_r(1,b)+C_{-r}(1,b)-b-1\over b+1}.
$$
Then
\begin{eqnarray}\disp\label{C1E03B}\nnb
 \dpp {} b F_0(r,b) &=&
  {rb^r+rb^{r-1}-b^r-1\over
(r-1)(b^r+1)(b+1)^2}  C_r(1,b)  \\
\disp \nnb & &+ {r +rb^{ -1}+b^r+1 \over
(r+1)(b^r+1)(b+1)^2}  C_{-r}(1,b)\\
\disp \nnb & & -{C_r(1,b)+C_{-r}(1,b) \over (b+1)^2}\\
\nnb\disp &=&  {r(b^{r-1}-1)C_{-r}(1,b)\over
 (r-1)(b^r+1)(b+1)^2} \Big[ {C_r(1,b)\over C_{-r}(1,b)} - {(r-1)(b^r-b^{ -1})\over  (r+1)(b^{r-1}-1)} \Big].
\end{eqnarray}
Consider
$$
F_1(r,b)=\ln {C_r(1,b)\over C_{-r}(1,b)}-\ln {(r-1)(b^r-b^{
-1})\over (r+1)(b^{r-1}-1)}, \qq\all r\geq 1, b>1.
$$
We have
\begin{eqnarray*}
 \disp \dpp {F_1(r,b)}b & = &  {1\over
 C_r(1,b)}\dpp {C_r(1,b)} b -{1\over
 C_{-r}(1,b)}\dpp {C_{_r}(1,b)} b \\
 \disp & & + {(r-1)b^{r-2}\over b^{r-1}-1} -{
rb^{r-1}+b^{-2}\over b^r-b^{-1}} \\
 \disp
   &=& {rb^r+rb^{r-1}-b^r-1\over (r-1)(b^r+1)(b+1)}  - {r +rb^{ -1}+b^r+1 \over  (r+1)(b^r+1)(b+1)}\\
 \disp & & + {(r-1)b^{r-2}\over b^{r-1}-1} -{
rb^{r-1}+b^{-2}\over b^r-b^{-1}} \\
 \disp
   &=& {F_2(r,b)\over (r-1)(r+1)b^2(b^r+1)(b+1)(b^{r-1}-1)(b^r-b^{-1})},
\end{eqnarray*}
where
\begin{eqnarray*}
\disp F_2(r,b) & =&
  -(r-1) b^{3r+1}+
 (r+1)b^{3r}+ r^2( r-1)b^{2r+2} \\
 \disp &&  +
  ( r+1) ( r^2-4r+1)b^{2r+1}-( r-1) ( r^2+4r+1)b^{2r}\\
 \disp &&  -r^2( r+1) b^{2r-1}  +  r^2( r+1) b^{r+2} +
 ( r-1)( r^2+4r+1) b^{r+1} \\
 \disp &&   -( r+1)( r^2-4r+1) b^r-r^2( r-1) b^{r-1}-(r+1)b+ (r-1)\\
  & =&  2e^{(3r+1)x}\Big(
  -(r-1) \sh (3r+1)x+
 (r+1)\sh (3r-1)x \\
 \disp &&  + r^2( r-1)\sh (r+3)x+
  ( r+1) ( r^2-4r+1)\sh (r+1)x\\
 \disp & & -( r-1) ( r^2+4r+1)\sh (r-1)x-r^2( r+1) \sh (r-3)x\Big)\\
 &\equiv & 2 e^{(3r+1)x} G_2(r,x)
\end{eqnarray*}
and $\disp x=\ln \sqrt b$.
Let
\begin{eqnarray*}\disp
G_3(r,x) &  \defeq & \Big(\dpp {} x-(r-1)\Big)G_2(r,x),  \\
G_4(r,x) &  \defeq &  \Big(\dpp {} x+(r-1)\Big)G_3(r,x)= \Big(\ddp 2 {} x-(r-1)^2\Big)G_2(r,x)\\
G_5(r,x) &  \defeq & \Big(\dpp {} x-(r+1)\Big)G_4(r,x),  \\
G_6(r,x) &  \defeq &  \Big(\dpp {} x+(r-1)\Big)G_5(r,x)= \Big(\ddp 2 {} x-(r+1)^2\Big)G_4(r,x)\\
\disp
G_7(r,x) &  \defeq & \Big(\dpp {} x-(r-3)\Big)G_6(r,x),  \\
G_8(r,x) &  \defeq &  \Big(\dpp {} x+(r-3)\Big)G_7(r,x)= \Big(\ddp 2 {} x-(r-3)^2\Big)G_6(r,x)\\
\disp
G_9(r,x) &  \defeq & \Big(\dpp {} x-(r+3)\Big)G_8(r,x),  \\
G_{10}(r,x) &  \defeq &  \Big(\dpp {} x+(r+3)\Big)G_9(r,x)=\Big(\ddp 2 {} x-(r+3)^2\Big)G_8(r,x),\\
& & \qq\qq \qq \qq\qq\qq \qq\qq\all x>0, r>1.
\end{eqnarray*}
Denote
$$
\pmatrix{\Gl_1\cr \Gl_2\cr \Gl_3\cr \Gl_4\cr \Gl_5\cr \Gl_6}=\pmatrix{r-1\cr r+1\cr r-3\cr r+3\cr 3r-1\cr 3r+1},
\GL =\pmatrix{ \Gl_1  &  &    &   &  &  \cr
              &  \Gl_2  &    &   &  &  \cr
              &   & \Gl_3    &   &  &  \cr
              &   &    &  \Gl_4  &  &  \cr
              &   &  &   &  \Gl_5   &  \cr
              &   &   &   &  &   \Gl_6}, X(x)=\pmatrix{\sh \Gl_1 x\cr\sh \Gl_2x\cr   \sh \Gl_3 x\cr \sh \Gl_4 x\cr \sh \Gl_5 x\cr\sh \Gl_6 x} 
$$
and define
$$
A_2=\pmatrix{ -( r-1) ( r^2+4r+1) \cr ( r+1) ( r^2-4r+1)\cr -r^2( r+1) \cr  r^2( r-1)\cr (r+1)\cr -(r-1)
    }
$$
and
$$
A_{2(k+1)}^\top=A_{2k}^\top(\GL^2-\Gl_k^2 I_6), \qq k=1,2,3,4,
$$
where $I_n$ denotes the $n\times n$ unit matrix. Then we have
$$
G_{2k}(r,x)=A_{2k}^\top X(x), \qq k=1,2,3,4,5.
$$
Further, we can get that
\begin{equation}\label{zero}
G_k(r,0)=\left\{\begin{array}{ll}\disp
0, & k=2, 3,4,6,7,8,\\
16r(r-1)(r+1)^2,& k=5,\\
4(r+1)^2(r-1)^2,& k=7, \\
-4480r^3(r-1)^2(r+1)^2, & k=9\end{array}\right.
\end{equation}
and
\begin{eqnarray*}\disp
&& {G_{10}(r,x) \over 1024 r^2(r-1)^2 (r+1)^2(2r-1)(2r+1)} \\
 &=& -  (r+2) \sh (3r+1)x +(r-2)\sh (3r-1)x \\
&<& 0, \qq\qq\all x>0, r>1.
\end{eqnarray*}

Now, we call a function $g$  poses {\textbf{Property
$(S)$} on $(\Ga,+\infty)$  if
$$
 \exists\, A\in
(\Ga,+\infty), \q\mbox{such that}\q \left\{\begin{array}{ll} g(x)>0 & \eqin (\Ga,A) \\ g(x)<0 & \eqin (A,+\infty).\end{array}\right.
$$
Noting that $G_k(r,+\infty)=-\infty$
and
$$
F_1(r,1)=1,\q F_1(r,+\infty)=-\infty,
$$
we get from \refeq{zero} that for fixed $r>1$,
$$
\begin{array}{c}\disp
G_{10}(r,x)<0,\qq\all x>0  \\
 \disp
  \Downarrow \\
 \disp
 G_{9}(r,x)<0,\qq\all x>0 \\
 \disp
 \Downarrow\\
 \disp
G_8(r,x)<0,\qq\all x>0 \\
 \disp
 \Downarrow\\
 \disp
G_7(r,x) \q\mbox{poses Property (S) on } (0,+\infty)\\
 \disp
 \Downarrow\\
 \disp
G_6(r,x) \q\mbox{poses Property (S) on } (0,+\infty)\\
 \disp
 \Downarrow\\
 \disp
\vdots\\
 \disp
 \Downarrow\\
 \disp
G_2(r,x) \q\mbox{poses Property (S) on } (0,+\infty)\\
 \disp
 \Downarrow\\
 \disp
F_2(r,b) \q\mbox{poses Property (S) on } (1,+\infty)\\
 \disp
 \Downarrow\\
 \disp
F_1(r,b) \q\mbox{poses Property (S) on } (1,+\infty).
\end{array}
$$
\if{
$$
G_{10}(r,x)<0,\qq\all x>0,
$$
$$
 \Big\Downarrow
$$
$$
G_{9}(r,x)<0,\qq\all x>0,
$$
$$
 \Big\Downarrow
$$
$$
G_8(r,x)<0,\qq\all x>0,
$$
$$
 \Big\Downarrow
$$
$$
G_7(r,x) \q\mbox{poses Property (S) on } (0,+\infty)
$$
$$
 \Big\Downarrow
$$
$$
G_6(r,x) \q\mbox{poses Property (S) on } (0,+\infty)$$
$$
 \Big\Downarrow
$$
$$
\vdots
$$
$$
 \Big\Downarrow
$$
$$
G_2(r,x) \q\mbox{poses Property (S) on } (0,+\infty)
$$
$$
 \Big\Downarrow
$$
$$
F_2(r,b) \q\mbox{poses Property (S) on } (1,+\infty)
$$
$$
 \Big\Downarrow
$$
$$
F_1(r,b) \q\mbox{poses Property (S) on } (1,+\infty).
$$
}\fi
Therefore, for some $b_0=b_0(r)\in (1,+\infty)$, $F_0(r,b)$ is increasing
strictly in $b\in (1,b_0)$ and decreasing  strictly in $b\in (b_0,
+\infty)$. Consequently,
\begin{eqnarray*}\disp
& & F_0(r,b) > \min\Big(F_0(r,1), F_0(r,+\infty)\Big)\\
\disp &=& \min(0,0)=0, \qq\all b>1.
\end{eqnarray*}
That is, \refeq{E305} holds for $r>1$.

For the case of $r=1$, we can prove similarly that
$$
F_1(1,b) \q\mbox{poses Property (S) on } (1,+\infty)
$$
and then get \refeq{E305}. We can also prove \refeq{E305} for $r=1$ in the following manner.

First, we can verify that $F_1(1,+\infty)=-\infty$.  On the other hand, along a subsequence $r\to 1^+$,
$b_0(r)$ tends to $\ell$ with $\ell=0,+\infty$ or a positive number.

If $\ell=0$, then by the continuity of $F_1$  we have $F_1(1,b)\leq 0$. Then $F_0(1,+\infty)<F_0(1,1)$ since $F_1(1,b)$ is negative for large $b$. This contradicts to
$F_0(1,+\infty)=F_0(1,1)=0$.

If $\ell=+\infty$, then  by the continuity of $F_1$  we have $F_1(1,b)\geq 0$.  This contradicts to $F_1(1,+\infty)=-\infty$.

Thus, we must have $\ell\in (0,+\infty)$ and
$$
\left\{\begin{array}{ll}\disp F_1(1,b)\geq 0, & \eqif b\in (0,\ell),\\
F_1(1,b)\leq 0, & \eqif b\in (\ell,+\infty).\end{array}\right.
$$
Therefore $F_0(1,b)$ is increasing
 in $b\in (1,\ell)$ and decreasing   in $b\in (\ell,
+\infty)$. Finally, we can get  \refeq{E305} since $F_0(1,b)$ is analytic and not a constant in $(1,+\infty)$.
\endpf

\def\theequation{3.\arabic{equation}}
\setcounter{equation}{0} \setcounter{section}{3}
\setcounter{Definition}{0} \setcounter{Example}{0}\textbf{3. Proof
of Conjecture A.}  We will prove Conjecture A in this section.
\if{First, we recall a basic property of $L_r(a,b)$ (see \cite{S}).
\begin{Proposition}\tt\label{T200} Let $b>a>0$, $-\infty<s<r<+\infty$. Then
\begin{equation}\label{E401}
\min (a,b)<L_s(a,b)<L_r(a,b)<\max (a,b).
\end{equation}
\end{Proposition}

Now we turns to prove Conjecture A.}\fi
\begin{Theorem}\tt\label{T402AA} Let $r\in (0,+\infty]$, $b>a>0$. Then
we have
\begin{equation}\label{E403A}
L_r(a,b)+L_{-r}(a,b)>2L_0(a,b).
\end{equation}
\end{Theorem}
\proof We need only to consider the cases of $r\in (0,+\infty)$ since
 \refeq{E403A} holds obviously when $r=+\infty$:
$$
L_{+\infty}(a,b)+L_{-\infty}(a,b)=2L_2(a,b)>2L_0(a,b).
$$
 Moreover, we can suppose that $b>a=1$ without loss of generality.

By \refeq{E110},  there exists a $\Gb=\Gb(r) >1$ such that
\begin{equation}\label{E404A}
L_r(1,b)+L_{-r}(1,b)>2L_0(1,b), \qq\all b\in (1,\Gb).
\end{equation}
Thus by \refeq{E104},Theorem \ref{T303},  Propositions \ref{T200} and  \ref{T301}, we
have that, for any $b\in (1,\Gb)$,
\begin{eqnarray*}\disp
&& L_r(1,b^2)+L_{-r}(1,b^2) \\
\disp &=& L_r(1,b)C_r(1,b)+L_{-r}(1,b)C_r(1,b)\\
\disp &=& {1\over
2}\Big(L_r(1,b)+L_{-r}(1,b)\Big)\Big(C_r(1,b)+C_{-r}(1,b)\Big)\\
\disp & & +{1\over
2}\Big(L_r(1,b)-L_{-r}(1,b)\Big)\Big(C_r(1,b)-C_{-r}(1,b)\Big)\\
\disp &>& {1\over
2}\Big(L_r(1,b)+L_{-r}(1,b)\Big)\Big(C_r(1,b)+C_{-r}(1,b)\Big)\\
\disp &>& 2L_0(1,b)C_0(1,b)=2L_0(1,b^2).
\end{eqnarray*}
Therefore,
\begin{equation}\label{E405A}
L_r(1,b)+L_{-r}(1,b)>2L_0(1,b), \qq\all b\in (1,\Gb^2).
\end{equation}
By induction, we can get that
\begin{equation}\label{E405BB}
L_r(1,b)+L_{-r}(1,b)>2L_0(1,b), \qq\all b\in (1,+\infty).
\end{equation}
We get the proof.
\endpf

\def\theequation{4.\arabic{equation}}
\setcounter{equation}{0} \setcounter{section}{4}
\setcounter{Definition}{0} \setcounter{Example}{0}\textbf{4.
 Proof
of Conjecture 2. } This section devotes to prove Conjecture 2. We state a lemma first.
\begin{Lemma}\tt\label{T401} Let $r\in (1,+\infty)$, $b>a>0$. Then
\begin{equation}\label{E401A}
 {b^2+a^2\over (b+a)^2}  < {b^{2r}+a^{2r}\over (b^r+a^r)^2}.
\end{equation}
Equivalently,
\begin{equation}\label{E401}
C^2_r(a,b)<C_r(a^2,b^2).
\end{equation}
\end{Lemma}
\proof The lemma follows directly from that $\disp {x^2+a^2\over (x+a)^2}$ is increasing strictly in $x\in (a,+\infty)$. \endpf

Now we state Conjecture 2 as
\begin{Theorem}\tt\label{T402} Let $r\in [0,+\infty)$, $b>a>0$. Then
\begin{equation}\label{E403}
C^2_r(a,b)+C^2_{-r}(a,b)<a^2+b^2.
\end{equation}
\end{Theorem}
\proof Without loss of generality, assume that $b>a=1$. We will prove \refeq{E403} by discussing four cases.

\textbf{Case I:  $r\in [1,3]$.}   As we pointed out in \cite{Lou},  by Proposition \ref{T301}(iii),
\begin{eqnarray*}
C^2_r(1,b)+C^2_{-r}(1,b) &<& C^2_3(1,b)+C^2_{-1}(1,b)= b^2+1.
\end{eqnarray*}

\textbf{Case II: $r\in (3,+\infty)$\footnote{In fact, this step holds for all $r\in (1,+\infty)$.}.}
Denote
$$
H_0(r,b)={C^2_r(1,b)+C^2_{-r}(1,b)-b^2-1\over 2(b^2+1)}.
$$
By \refeq{E306}---\refeq{E307},
\begin{eqnarray}\disp\label{C1E03}\nnb
 && \dpp {} b H_0(r,b) \\
\nnb &=&
   {rb^r+rb^{r-1}-b^r-1\over
(r-1)(b^r+1)(b+1)(b^2+1)} C^2_r(1,b)  \\
\disp \nnb & &+ {r +rb^{ -1}+b^r+1 \over
(r+1)(b^r+1)(b+1)(b^2+1)}C^2_{-r}(1,b)  \\
\disp \nnb & & -{b\over  (b^2+1)^2}(C^2_r(1,b)+C^2_{-r}(1,b))\\
\nnb\disp &=&
   {b^{r+1}+(r-1)b^r+rb^{r-1}-rb^2-(r-1)b-1\over
(r-1)(b+1)(b^2+1)^2(b^r+1)} C^2_r(1,b)  \\
\disp \nnb & &- {rb^{r+3} +(r+1)b^{r+2}-b^{r+1}+b^2-(r+1)b-r \over
(r+1)b(b+1)(b^2+1)^2(b^r+1)}C^2_{-r}(1,b)  \\
\nnb\disp &=&
   {b^{r+1}+(r-1)b^r+rb^{r-1}-rb^2-(r-1)b-1\over
(r-1)(b+1)(b^2+1)^2(b^r+1)} \Big({C^2_r(1,b)\over C^2_{-r}(1,b)}  \\
\disp & &- {(r-1)  (rb^{r+2} +(r+1)b^{r+1}-b^r+b-(r+1)-rb^{-1}) \over
(r+1) (b^{r+1}+(r-1)b^r+rb^{r-1}-rb^2-(r-1)b-1)} \Big).
\end{eqnarray}
Obviously, for any $r\geq 3$, $b>1$, it holds that\footnote{It is not very hard but a little complex to prove that \refeq{E406} holds for all $r>1, b>1$.}
\begin{eqnarray}\label{E406}
  && b^{r+1}+(r-1)b^r+rb^{r-1}-rb^2-(r-1)b-1>0.
\end{eqnarray}
\if{Generally, we have
\begin{eqnarray}\label{E407}
\nnb && b^{r+1}+(r-1)b^{r}+rb^{r-1}-rb^{2}-(r-1)b-1\\
\nnb &=& b^2\int_1^b dx\int_1^x dy\int_1^y 2r(r-1)x^{r-4}(x-y)y^{-r}z^{-2}(z^{r+1}-1)\, dz\\
\nnb && +(r-1)(r^2-r+4)b^2\int_1^b dx
\int_1^x x^{r-4}(x-y)y^{-r}\, dy+b^2\int_1^b
 2(r-1)^2x^{r-4}\, dx\\
&>& 0, \qq\all r>1, b>1.
\end{eqnarray}}\fi
On the other hand,
\begin{eqnarray}\label{E408}
\nnb && rb^{r+2}+(r+1)b^{r+1}-b^{r}+b-(r+1)-rb^{-1}\\
\nnb &=& b^r(b-1)+r(b^{r+2}-b^{-1})+rb^{r+1}+b-(r+1)\\
&>& 0, \qq\all r>1, b>1.
\end{eqnarray}
Thus we can define
\begin{eqnarray}\label{E408B}
\nnb &&
H_1(r,b)=\ln {C^2_r(1,b)\over C^2_{-r}(1,b)} \\
&&  -\ln  {(r-1) (rb^{r+2} +(r+1)b^{r+1}-b^r+b-(r+1)-rb^{-1})  \over
(r+1) (b^{r+1}+(r-1)b^r+rb^{r-1}-rb^2-(r-1)b-1)}.
\end{eqnarray}
We have
\begin{eqnarray}\label{EH1}
\nnb\dpp {} b H_1(r,b)&=& 2{rb^r+rb^{r-1}-b^r-1\over
(r-1)(b^r+1)(b+1)} \\
\nnb && - 2{r +rb^{ -1}+b^r+1 \over
(r+1)(b^r+1)(b+1)} \\
\nnb & &  -   {  r(r+2)b^{r+1} +(r+1)^2b^r-rb^{r-1}+1+rb^{-2}  \over   rb^{r+2} +(r+1)b^{r+1}-b^r+b-(r+1)-rb^{-1} } \\
\nnb & & + {(r+1)b^r+(r-1)rb^{r-1}+r(r-1)b^{r-2}-2rb-(r-1)  \over  b^{r+1}+(r-1)b^r+rb^{r-1}-rb^2-(r-1)b-1 }\\
\nnb &=& {1\over b^{r+1}+(r-1)b^r+rb^{r-1}-rb^2-(r-1)b-1} \\
\nnb & & \cdot {1\over  rb^{r+2} +(r+1)b^{r+1}-b^r+b-(r+1)-rb^{-1}} \\
& & \cdot  {rH_2(r,b)\over (r-1)(r+1)b^2(b+1)(b^r+1)},
\end{eqnarray}
where
\begin{eqnarray}\label{EW2}
\nnb H_2(r,b)&=&  (r-1)^2b^{3r+5}+(r^2+6r-3)b^{3r+4}+4(r+1)b^{3r+3}+4(r-1)b^{3r+2}\\
\nnb & & -(r^2-6r-3)b^{3r+1}-(r+1)^2b^{3r}+r^2(r-1)^2b^{2r+6}\\
\nnb & & +(3r^4-6r^3-6r^2-2r+3)b^{2r+5} +(3r^4-10r^3-6r^2+6r-9)b^{2r+4}\\
\nnb  & & +(r^4-14r^3+r^2+4r+12)b^{2r+3} -(r^4+14r^3+r^2-4r+12)b^{2r+2}\\
\nnb & & -(3r^4+10r^3-6r^2-6r-9)b^{2r+1} -(3r^4+6r^3-6r^2+2r+3)b^{2r}\\
\nnb & & -r^2(r+1)^2b^{2r-1} +r^2(r+1)^2b^{r+6}+(3r^4+6r^3-6r^2+2r+3)b^{r+5}\\
\nnb && +(3r^4+10r^3-6r^2-6r-9)b^{r+4}+(r^4+14r^3+r^2-4r+12)b^{r+3}  \\
\nnb && -(r^4-14r^3+r^2+4r+12)b^{r+2}-(3r^4-10r^3-6r^2+6r-9)b^{r+1}\\
\nnb  & & -(3r^4-6r^3-6r^2-2r+3)b^r-r^2(r-1)^2b^{r-1}+(r+1)^2b^5 \\
\nnb   & & +(r^2-6r-3)b^4-4(r-1)b^3-4(r+1)b^2-(r^2+6r-3)b\\
\nnb   & & -(r-1)^2\\
\nnb &=& 2e^{(3r+5)x} \Big((r-1)^2\sh (3r+5)x+(r^2+6r-3)\sh (3r+3)x \\
\nnb && +4(r+1)\sh (3r+1)x +4(r-1)\sh (3r-1)x\\
\nnb & & -(r^2-6r-3)\sh (3r-3)x-(r+1)^2\sh (3r-5)x\\
\nnb & & +r^2(r-1)^2\sh (r+7)x  +(3r^4-6r^3-6r^2-2r+3)\sh (r+5)x \\
\nnb & & +(3r^4-10r^3-6r^2+6r-9)\sh (r+3)x\\
\nnb & &  +(r^4-14r^3+r^2+4r+12)\sh (r+1)x\\
\nnb & &  -(r^4+14r^3+r^2-4r+12)\sh (r-1)x \\
\nnb & &  -(3r^4+10r^3-6r^2-6r-9)\sh (r-3)x \\
\nnb & & -(3r^4+6r^3-6r^2+2r+3)\sh (r-5)x  -r^2(r+1)^2\sh (r-7)x\Big)\\
 &\equiv  & 2e^{(3r+5)x} W_2(r,x),
\end{eqnarray}
where $x=\ln \sqrt b$.  Similar to Section 2, we define
\if{$$
\mu=\pmatrix{\mu_1\cr\mu_2\cr\mu_3\cr \mu_4\cr \mu_5\cr \mu_6\cr \mu_7\cr \mu_8\cr \mu_9\cr \mu_{10}\cr
\mu_{11}\cr \mu_{12}\cr \mu_{13}\cr \mu_{14}}=\pmatrix{r-1\cr r+1\cr r-3\cr r+3\cr r-5\cr r+5\cr 3r-1\cr 3r+1\cr r-7\cr 3r-3\cr 3r+3\cr r+7\cr 3r-5\cr 3r+5},
$$
\begin{eqnarray*}\disp
W_{2k+1}(r,x) &  \defeq & \Big(\dpp {} x-\mu_k \Big)W_{2k}(r,x), \q k=1,2,\ldots, 13, \\
W_{2k+2}(r,x) &  \defeq & \Big(\dpp {} x+\mu_k\Big)W_{2k+1}(r,x), 
\q k=1,2,\ldots, 13.
\end{eqnarray*}}\fi
$$
\pmatrix{\mu_1  & \mu_8\cr \mu_2 & \mu_9 \cr\mu_3 & \mu_{10} \cr \mu_4 & \mu_{11}\cr \mu_5 & \mu_{12}\cr \mu_6 & \mu_{13}\cr \mu_7 & \mu_{14} }=\pmatrix{r-1 & 3r+1\cr r+1 & r-7\cr r-3 & 3r-3 \cr r+3& 3r+3\cr r-5& r+7\cr r+5& 3r-5\cr 3r-1&  3r+5},
$$
\begin{eqnarray*}\disp
\left\{\begin{array}{l} \disp W_{2k+1}(r,x)  \defeq  \Big(\dpp {} x-\mu_k \Big)W_{2k}(r,x), \vspace{2mm}\\
\disp W_{2k+2}(r,x)   \defeq  \Big(\dpp {} x+\mu_k\Big)W_{2k+1}(r,x),  \end{array}\right.
\q k=1,2,\ldots, 13.
\end{eqnarray*}
\if{
\begin{eqnarray*}\disp
W_3(r,x) &  \defeq & \Big(\dpp {} x-(r-1)\Big)W_2(r,x),  \\
W_4(r,x) &  \defeq & \Big(\dpp {} x+(r-1)\Big)W_3(r,x)=  \Big(\ddp 2 {} x-(r-1)^2\Big)W_2(r,x)\\
W_5(r,x) &  \defeq & \Big(\dpp {} x-(r+1)\Big)W_4(r,x),  \\
W_6(r,x) &  \defeq &  \Big(\dpp {} x+(r+1)\Big)W_5(r,x)=  \Big(\ddp 2 {} x-(r+1)^2\Big)W_4(r,x)\\
W_7(r,x) &  \defeq & \Big(\dpp {} x-(r-3)\Big)W_6(r,x),  \\
W_8(r,x) &  \defeq &  \Big(\dpp {} x+(r-3)\Big)W_7(r,x)=  \Big(\ddp 2 {} x-(r-3)^2\Big)W_6(r,x)\\
W_9(r,x) &  \defeq & \Big(\dpp {} x-(r+3)\Big)W_8(r,x),  \\
W_{10}(r,x) &  \defeq &  \Big(\dpp {} x+(r+3)\Big)W_9(r,x)=  \Big(\ddp 2 {} x-(r+3)^2\Big)W_8(r,x)\\
W_{11}(r,x) &  \defeq & \Big(\dpp {} x-(r-5)\Big)W_{10}(r,x),  \\
W_{12}(r,x) &  \defeq &  \Big(\dpp {} x+(r-5)\Big)W_{11}(r,x)=  \Big(\ddp 2 {} x-(r-5)^2\Big)W_{10}(r,x)\\
\end{eqnarray*}

\begin{equation}\label{zero}
W_k(r,0)=\left\{\begin{array}{ll}\disp
0, & k=2, 3,4,6,7,8,10,12,14,\\
128 (r-1)^2(r+1)^2(r^2+3), & k=5,\\
18432 (r-1)^2(r+1)^2(r^2+1),& k=7, \\
2048 (r-1)^2(r+1)^2(49r^4+54r^3+699r^2+54r+180), & k=9,\\
32768 (r-1)^2(r+1)^2 (26r^6+535r^4+2019r^2+180), & k=11,\\
16384 (r-1)^2(r+1)^2 (488r^7+520r^6+14131r^5+10700r^4+63873r^3+40380r^2+94080r+3600), & k=13,\\
\end{array}\right.
\end{equation} }\fi
We have
\begin{eqnarray*}
W_{2k}(r,0)& =& 0, \q k=1, 2,\ldots, 12,\\
W_3(r,0)& =& 0, \\
W_5(r,0)& =& 128 (r-1)^2(r+1)^2(r^2+3),  \\
W_5(r,0)& =& 128 (r-1)^2(r+1)^2(r^2+3),  \\
W_7(r,0)& =&  18432 (r-1)^2(r+1)^2(r^2+1),   \\
W_9(r,0)& =& 2048 (r-1)^2(r+1)^2(49r^4+54r^3+699r^2+54r+180),  \\
W_{11}(r,0)& =& 32768 (r-1)^2(r+1)^2 (26r^6+535r^4+2019r^2+180),   \\
W_{13}(r,0)& =& 16384 (r-1)^2(r+1)^2 (488r^7+520r^6 +14131r^5  \\
& &  +10700r^4+63873r^3+40380r^2+94080r+3600),   \\
W_{15}(r,0)& =& 262144r^{2}(r-1)^2(r+1)^2 (280r^8  +11536r^6 +67429r^4  \\
& & +103185r^2+117612),   \\
W_{17}(r,0)& =& 1572864 r^2(r-1)^2(r+1)^2
(48r^{10}280r^{9}+10246r^{8}\\
&& +11536r^{7}+181169r^{6} +67429r^{5}+584375r^{4}\\
& & +103185r^{3}+725272r^{2}+117612r+924768),\\
W_{19}(r,0)& =& 25165824 r^2(r-1)^2(r+1)^2
(1674r^{10}+98095r^{8}\\
& & +912478r^{6}+1651021r^{4}+2096148r^{2}+2748384),   \\
W_{21}(r,0)& =& 50331648 r^2(r-1)^2(r+1)^2
(7668r^{12}+11718r^{11}\\
&& +600358r^{10}+686665r^{9}+7482173r^{8}+6387346r^{7}\\
& & +18275708r^{6}+11557147r^{5}+977689r^{4}+14673036r^{3}\\
& & +16973484r^{2}+19238688r-319680),\\
W_{23}(r,0)& =& 1207959552r^2(r-1)^2(r+1)^2(324r^{14}+6318r^{13}+115858r^{12}\\
& & +591119r^{11}+3938675r^{10}+8991557r^{9}+27578739r^{8}\\
& &+31160805r^{7} +37714913r^{6}+15674485r^{5}+1539103r^{4}\\
& & +40673476r^{3}+56552388r^{2}+31638240r-216000),\\
W_{25}(r,0)& =& 2415919104r^2(r-1)^2(r+1)^2(75956r^{14}
+118328r^{13}\\
&& +7598119r^{12}+8628046r^{11}+131265979r^{10}+114185470r^{9}\\
&& +559333873r^{8}+414634234r^{7}+392275685r^{6}+21748202r^{5}\\
& & -439324492r^{4}+596357720r^{3}+1115006880r^{2}+636048000r\\
& & -1728000),\\
W_{27}(r,0)& =& 38654705664 r^2(r-1)^3(r+1)^3(2r-1)(2r+1)(10692r^{12}\\
& & +1371339r^{10}+31347410r^{8}+183116951r^{6}+282237368r^{4} \\
& & +92029680r^{2}+2592000).
\end{eqnarray*}
While
\begin{eqnarray}\label{E404}\disp
\nnb W_{28}(r,x) &=& (r-1)^2\prod^{13}_{k=1} \Big((3r+5)^2-\mu_k^2\Big)\sh (3r+5)x\\
\nnb  &=& (r-1)^2\prod^{13}_{k=1} \Big((3r+5)^2-\mu_k^2\Big)\sh (3r+5)x \\
&>& 0, \qq\qq\all x>0, r>1
\end{eqnarray}
since $3r+5>|\mu_k|$  for $r>1$ and $1\leq k\leq 13$. It is easy to see from above that when $r>1$,
\begin{equation}\label{E405}
W_k(r,0)=\left\{\begin{array}{ll}\disp
0, & k=2m, \q m=1,2,\ldots, 13,\\
0, & k=3,\\
>0, & k=2m+1,\q m=2,3,\ldots, 13.
\end{array}\right.
\end{equation}
Combining \refeq{E405} with \refeq{E404} we get that
\begin{equation}\label{E405B}
 W_2(r,x)>0, \qq\all r>1,x>0.
\end{equation}
Combing \refeq{EH1}---\refeq{EW2} with \refeq{E405B}, we get that $H_1(r,b)$ is increasing strictly in $b\in [1,+\infty)$.
Thus there exists a $b_1=b_1(r)\in (1,+\infty)$ such that
$H_1(r,b)$ is negative in $(1,b_1)$ and positive in $(b_1,+\infty)$ since
$$
H_1(r,1)=-\ln {r+1\over r-1}<0,\q H_1(r,+\infty)=+\infty.
$$
Consequently, $H_0(r,b)$ is decreasing strictly in $(1,b_1)$ and increasing strictly in $(b_1,+\infty)$.
Therefore
$$
H_0(r,b)<\min(H_0(r,1),H_0(r,+\infty))=\min(0,0)=0, \qq\all r>1,b>1.
$$
That is,
$$
C^2_r(1,b)+C^2_{-r}(1,b) <1+ b^2,\qq \all r>1,b>1.
$$

\textbf{Case III: $r=0$.} We have
$$
2C^2_0(1,b)=2\Big({1+b\over 2}\Big)^2 <1+ b^2,\qq \all b>1.
$$

\textbf{Case IV: $r\in (0,1)$.}  Denote $\disp s={1\over r}$. Then $s>1$. By Lemma \ref{T401} and what we got in Case II, we have
\begin{eqnarray*}
&& C^2_r(1,b^s)+C^2_{-r}(1,b^s) \\
&=&  \Big({1+b^s\over 1+b}\Big)^2C^2_s(1,b)+\Big({1+b^s\over 1+b}\Big)^2 C^2_{-s}(1,b)\\
&<& \Big({1+b^s\over 1+b}\Big)^2(1+b^2) < 1+b^{2s}, \qq\all b>1.
\end{eqnarray*}
Therefore,
$$
 C^2_r(1,b)+C^2_{-r}(1,b)<1+b^2, \qq\all b>1.
$$

Combining Cases I---IV, we get the proof.
\endpf

One can get immediately from Theorem \ref{T402} that
\begin{Corollary}\tt\label{T402B} Let $r\in [0,+\infty)$, $a,b>0$, $a\ne b$. Then
we have
\begin{equation}\label{E410}
C_{+\infty}(a,b)C_r(a,b)+C_{-\infty}(a,b)C_{-r}(a,b)<a^2+b^2.
\end{equation}
\end{Corollary}
\if{\proof For $b>a=1$, using  we have
\begin{eqnarray*}
& & bC_r(1,b)+ C_{-r}(1,b)\\
&\leq & {1\over 2}\Big(b^2+C^2_r(1,b)\Big)+{1\over 2}\Big(1+C^2_{-r}(1,b)\Big)\\
&<&  b^2+1.
\end{eqnarray*}
Then \refeq{E410} follows.
\endpf}\fi

\def\theequation{5.\arabic{equation}}
\setcounter{equation}{0} \setcounter{section}{5}
\setcounter{Definition}{0} \setcounter{Example}{0}\textbf{5. Proof
of Conjecture B.}

We turn to prove Conjecture B and state it as
\begin{Theorem}\tt\label{T601A} Let $r\in [0,+\infty)$, $b>a>0$. Then
we have
\begin{equation}\label{E601A}
L_r(a,b)+L_{-r}(a,b)<a+b.
\end{equation}
\end{Theorem}
\proof Without loss of generality, we suppose that $b>a=1$. Since
\refeq{E601A} holds obviously for $r=0$, we suppose that $r\in (0,+\infty)$ in the following.
By \refeq{E110A}, there exists a $\Gg>1$, such that
\begin{equation}\label{E601AA}
L_r(1,b)+L_{-r}(1,b)<b+1, \qq\all b\in (1,\Gg).
\end{equation}
Thus by Corollary  \ref{T402}, Propositions \ref{T200}--\ref{T301}, and noting that
$$
L_r(1,b)-L_{-r}(1,b)<b-1,
$$
we have
\begin{eqnarray*}\disp
&& L_r(1,b^2)+L_{-r}(1,b^2) \\
\disp &=& L_r(1,b)C_r(1,b)+L_{-r}(1,b)C_r(1,b)\\
\disp &=& {1\over 2}\Big(L_r(1,b)+L_{-r}(1,b)\Big)\Big(C_r(1,b)+C_{-r}(1,b)\Big)\\
\disp & & +{1\over 2}\Big(L_r(1,b)-L_{-r}(1,b)\Big)\Big(C_r(1,b)-C_{-r}(1,b)\Big)\\
\disp &<& {1\over 2} (b+1)\Big(C_r(1,b)+C_{-r}(1,b)\Big)\\
\disp & & +{1\over 2}(b-1)\Big(C_r(1,b)-C_{-r}(1,b)\Big)\\
&=& bC_r(1,b)+C_{-r}(1,b)\\
\disp &<& b^2+1.
\end{eqnarray*}
Therefore,
\begin{equation}\label{E503}
L_r(1,b)+L_{-r}(1,b)<b+1, \qq\all b\in (1,\Gg^2).
\end{equation}
We get the proof by induction. \endpf

\def\theequation{6.\arabic{equation}}
\setcounter{equation}{0} \setcounter{section}{6}
\setcounter{Definition}{0}\setcounter{Example}{0} \textbf{6.
Further Results.} In this section, we will yield some related results.  \if{We recall a result proved by Lou in \cite{Lou}:
\begin{Lemma}\tt\label{T601} Let  $0<r<s<+\infty$, $b>a>0$. Then
\begin{equation}\label{E601}
ab<C_s(a,b)C_{-s}(a,b)<C_r(a,b)C_{-r}(a,b)< \Big({a+b\over 2}\Big)^2.
\end{equation}
\end{Lemma}
We mention that
$$
ab=C_{+\infty}(a,b)C_{-\infty}(a,b)=C^2_{-1}(a,b).
$$
}\fi
We have
\begin{Corollary}\tt\label{T603} Let  $0<r<+\infty$, $b>a>0$. Then
\begin{equation}\label{E604}
C_r(a,b)+C_{-r}(a,b)<\sqrt{3a^2+2ab+3b^3\over 2}
\end{equation}
and
\begin{equation}\label{E605}
C^2_r(a,b)+C^2_{-r}(a,b)>{(a+b)^2\over 2}.
\end{equation}
\end{Corollary}
\proof  We have
\begin{eqnarray*}
&& C_r(a,b)+C_{-r}(a,b) =\sqrt{C^2_r(a,b)+C^2_{-r}(a,b)+2C_r(a,b)C_{-r}(a,b)}\\
&<& \sqrt{a^2+b^2+2\Big({a+b\over 2}\Big)^2}=\sqrt{{3a^2+2ab+3b^3\over 2}}
\end{eqnarray*}
and
\begin{eqnarray*}
&& C^2_r(a,b)+C^2_{-r}(a,b)\geq {(C_r(a,b)+C_{-r}(a,b))^2\over 2}> {(a+b)^2\over 2}.
\end{eqnarray*}
We get the proof. \endpf

\begin{Remark}\tt By  Proposition \ref{T301}, Theorem \ref{T303} and Corollary \ref{T603}, for $0<r<+\infty$, $b>a>0$, we have the following inequalities:
\begin{eqnarray}\label{E605A}
\nnb & & C_{-1}(a,b)<\sqrt{C_r(a,b)C_{-r}(a,b)}<C_0(a,b)<{C_r(a,b)+C_{-r}(a,b)\over 2} \\
 &<& \sqrt{{C_0(a,b)(C_0(a,b)+C_2(a,b))\over 2}}.
\end{eqnarray}
\end{Remark}

The following result can be looked a corollary of Proposition \ref{T301}:
\begin{Lemma}\tt\label{T602} Let  $0<r<s<+\infty$, $b>a>0$. Then
\begin{equation}\label{E601B}
ab<L_s(a,b)L_{-s}(a,b)<L_r(a,b)L_{-r}(a,b)< L_0^2(a,b).
\end{equation}
\end{Lemma}
\proof Let $0<r<s<+\infty$, $b>a=1$. We have
$$
\lim_{b\to 1}{L_r(1,b)L_{-r}(1,b)-L_s(1,b)L_{-s}(1,b)\over (b-1)^4}={s^2-r^2\over 1440}.
$$
Thus, there exists a $\mu=\mu(r,s)>1$ such that for any $ b\in (1,\mu)$,
\begin{equation}\label{E602}
L_r(1,b)L_{-r}(1,b)>L_s(1,b)L_{-s}(1,b).
\end{equation}
Consequently,
\begin{eqnarray*}
&& L_r(1,b^2)L_{-r}(1,b^2)=L_r(1,b)L_{-r}(1,b)C_r(1,b)C_{-r}(1,b)\\
&>& L_s(1,b)L_{-s}(1,b)C_s(1,b)C_{-s}(1,b)=L_s(1,b^2)L_{-s}(1,b^2).
\end{eqnarray*}
That is, \refeq{E602} holds for $b\in (1,\mu^2)$. Thus, by induction, \refeq{E602} holds for $b\in (1,+\infty)$.
Moreover, it follows from \refeq{E602} that
\begin{eqnarray*}
&& L_s(1,b^2)L_{-s}(1,b^2)>\lim_{t\to +\infty} L_t(1,b)L_{-t}(1,b)=b
\end{eqnarray*}
and
\begin{eqnarray*}
&& L_r(1,b^2)L_{-r}(1,b^2)<\lim_{t\to 0} L_t(1,b)L_{-t}(1,b)=L_0^2(1,b).
\end{eqnarray*}
We get the proof. \endpf

On the other hand, we have:
\begin{Corollary}\tt\label{T604} Let  $0<r<+\infty$, $b>a>0$. Then
\begin{enumerate}\renewcommand{\labelenumi}{{\rm (\roman{enumi})}}
\item for any $\Ga\in (0,1]$,
\begin{equation}\label{E606}
C^\Ga_r(a,b)+C^\Ga_{-r}(a,b)>a^\Ga+b^\Ga;
\end{equation}
\item for any $\Gb\in [0,+\infty)$,
\begin{equation}\label{E607}
C^{2+\Gb}_r(a,b)+C^{2+\Gb}_{-r}(a,b)<a^{2+\Gb}+b^{2+\Gb};
\end{equation}
\item for any $\Gb\in [0,+\infty)$,
\begin{equation}\label{E608}
L^{1+\Gb}_r(a,b)+L^{1+\Gb}_{-r}(a,b)<a^{1+\Gb}+b^{1+\Gb}.
\end{equation}
\end{enumerate}
\end{Corollary}
\proof
For any $\Gb\in (0,+\infty)$, $b>a>0$, it is easy to prove that
\begin{eqnarray}\label{E609}
\nnb && x^{1+\Gb}+y^{1+\Gb}< a^{1+\Gb}+b^{1+\Gb},\\
& & \qq\qq\all (x,y)\in \{(u,v)| \sqrt{ab}<\sqrt{uv}\leq {u+v\over 2}<{a+b\over 2}\}.
\end{eqnarray}
\begin{enumerate}\renewcommand{\labelenumi}{{\rm (\roman{enumi})}}
\item Let  $0<r<+\infty$, $b>a>0$, $\Ga\in (0,1]$. We claim \refeq{E606} holds. Otherwise,
\begin{equation}\label{E611}
C^\Ga_r(a,b)+C^\Ga_{-r}(a,b)\leq a^\Ga+b^\Ga.
\end{equation}
Thus, by \refeq{E304} and take limitation in \refeq{E609}, we get
$$
C_r(a,b)+C_{-r}(a,b)\leq a+b.
$$
Contradicts to \refeq{E305}. Therefore, \refeq{E606} holds.
\item By  \refeq{E304} and \refeq{E403},
$$
C^2_r(a,b)C^2_{-r}(a,b)>a^2b^2, \q C^2_r(a,b)+C^2_{-r}(a,b)<a^2+b^2.
$$
Thus, it follows from \refeq{E609} that
$$
\Big(C^2_r(a,b)\Big)^{1+{\Gb\over 2}}+\Big(C^2_{-r}(a,b)\Big)^{1+{\Gb\over 2}}<(a^2)^{1+{\Gb\over 2}}+(b^2)^{1+{\Gb\over 2}}.
$$
That is, \refeq{E607} holds.
\item Similar to \refeq{E607}, we can get \refeq{E608} directly from \refeq{E103A}, \refeq{E601A}  and \refeq{E609}. \endpf
\end{enumerate}

\footnotesize
\vspace{5mm}\footnotesize \ \\

\end{document}